\documentclass{amsart}
\usepackage[pagewise]{lineno}
\usepackage{enumerate}
\usepackage{amsmath}
\usepackage{amssymb}
\usepackage{graphicx}
\usepackage{float}
\usepackage{mathrsfs}
\usepackage{caption}

\newtheorem{theorem}{Theorem}[section]

\newtheorem{proposition}[theorem]{Proposition}

\theoremstyle{definition}
\newtheorem{definition}[theorem]{Definition}
\newtheorem{example}[theorem]{Example}
\newtheorem{corollary}[theorem]{Corollary}

\theoremstyle{remark}
\newtheorem{remark}[theorem]{Remark}

\numberwithin{equation}{section}

\begin{document}

\title[The quantale of order-preserving maps]
{THE QUANTALE OF ORDER-PRESERVING MAPS }

\author{Hongwei Wu}
\address {School of Mathematics and Statistics, Shaanxi Normal University, Xi'an, 710119, P.R. China }
\email{1269991632@qq.com, wuhw@snnu.edu.cn}






\begin{abstract}
In this paper,  two new composition operations  are defined among the order-preserving maps. They can act on order-preserving maps like the usual composition operation.  They are coincide with the usual composition operation when  the order-preserving maps are sup-preserving maps or meet-preserving maps. The usual composition operation can't endow the set $L^{L}$ of order preserving maps on a complete lattices $L$ with quantale or co-quantale structures. Luckily, two new operations can endow $L^{L}$ with a quantale and a co-quantale structures respectively, whenever $L$ is a completely distributive lattice.
\end{abstract}

\maketitle



\section{Introduction}
\quad Quantales were introduced by Mulvey (see \cite{M86})  in order to provide a lattice-theoretic setting for studying non-commutative $C^{\ast}$-algebras, as well as constructive foundations for quantum mechanics (see \cite{BRV89,KPRR89,KR04,MP01}). Quantales form an
important class of ordered algebraic structures, and interest in quantales was
stimulated by the fact (see \cite{Yetter1990}) that Girard quantales provide a sound and
complete class of models for linear intuitionistic logic, just as complete Heyting algebras model intuitionistic logic. More various types and aspects of quantales, which are not explained here, please refer to \cite{Rosenthal1990}.

Given a complete lattice $L$, the set $\mathcal{S}(L,L)$ of all sup-preserving maps and $\mathcal{M}(L,L)$ of all meet-preserving maps are complete lattice with the pointwise order. It is well known that $(\mathcal{S}(L,L),\circ)$ is a quantale, where $\circ$ is the composition of maps. It is an important quantale in quanatle theory.  In 1996, Paseka introduced the concept of simple quantale which is closely related to $C^{\ast}$-algebra (see \cite{Paseka1996}). From \cite{PK2000}, $(\mathcal{S}(L,L),\circ)$ is a widely used simple quantale. From \cite{KP08,PK2000}, $(\mathcal{S}(L,L),\circ)$ is a unital quantale and each quantale can be embedded into a unital quantale. Besides, it is also a Girard quantale whenever $L$ is a completely distributive lattice (see \cite{KP08}). Base on these results, $(\mathcal{S}(L,L),\circ)$ may play an important role in the study of the relationship between linear intuitionistic logic and $C^{\ast}$-algebra in the future. Let $L^{L}$ denote the set of all order-preserving maps on $L$. $L^{L}$ is a complete lattice with the pointwise order.  $\mathcal{S}(L,L)$  and  $\mathcal{M}(L,L)$ are sublattices of $L^{L}$.  $(\mathcal{S}(L,L),\circ)$ can form a quantale.  However,  $(L^{L},\circ)$ is not a quantale. In the section 3,  we introduce a new composition $\cdot$, which can work in the same way as $\circ$ usually works with order-preserving maps.
It can endow $L^{L}$  with a quantale structure, whenever $L$ is a completely distributive lattice. We prove that $(\mathcal{S}(L,L),\circ)$ is isomorphic to a quotient of $(L^{L},\cdot)$. In fact, Picado do a similar work on order-reversing maps (see \cite{Picado2005}). Co-quantale is a contravariant notion of quantale. It was introduced in \cite{Flagg97} to deal with topological spaces as pseudo metric spaces where the the distance take valued in this suitable quantale.
It was also used in valued logic (see \cite{RZ2021}). We prove that $(\mathcal{M}(L,L),\circ)$ is a
co-quantale. However,  $(L^{L},\circ)$ is not a co-quantale. In the section 4,  we introduce another new composition $\bullet$, which can work in the same way as $\circ$ usually works with order-preserving maps. It can endow $L^{L}$  with a co-quantale structure, whenever $L$ is a completely distributive lattice. Moreover, we prove that $(\mathcal{M}(L,L),\circ)$ is a quotient of $(L^{L},\bullet)$.

\section{Preliminaries }

We refer  to \cite{DP02} for lattice theory and to \cite{Rosenthal1990} for quantale theory. In this paper, we use $1$
to denote the top element and $0$ to denote the bottom element in a complete lattice.

Let $A$ and $B$ be complete lattices. The set of all order-preserving maps from $A$ to $B$ is denoted by $B^{A}$. $B^{A}$ is ordered pointwisely:
\[f\leq g\ \ \ \mbox{in}\ B^{A}\Longleftrightarrow f(a)\leq g(a)\ \ \ \mbox{for\ all}\ a\in A.\]
$B^{A}$ is a complete lattice with the above order relation. The set of all sup-preserving maps from $A$ to $B$ is denoted by $\mathcal{S}(A,B)$ and the set of all meet-preserving maps from $A$ to $B$ is denoted by $\mathcal{M}(A,B)$. For  simplicity, we use $\mathcal{S}(A)$ to denote $\mathcal{S}(A,A)$ and use $\mathcal{M}(A)$ to denote $\mathcal{M}(A,A)$. $\mathcal{S}(A,B)$ and $\mathcal{M}(A,B)$ are all complete lattices with pointwise order. The symbol $\circ$ denotes the composition of maps.
\begin{remark}\label{rem2.1} (1) For all $\{f_{i}\}_{i\in I}\subseteq \mathcal{S}(A,B)$,\ $(\bigvee f_{i})(a)=\bigvee_{i\in I}f_{i}(a)$.

(2) For all $\{f_{i}\}_{i\in I}\subseteq \mathcal{M}(A,B)$,\ $(\bigwedge f_{i})(a)=\bigwedge_{i\in I}f_{i}(a)$.

(3) $f\in \mathcal{S}(A,B)$ and $g\in \mathcal{S}(B,C)$, then $g\circ f\in \mathcal{S}(A,C)$.

(4) Define a map $\underline{\perp}\colon A\longrightarrow B$ by $\underline{\perp}(a)=0_{B}$ for all $a\in A$. Then $\underline{\perp}$ is the bottom element in $B^{A}$. Define a map $\overline{\top}\colon A\longrightarrow B$ by $\overline{\top}(a)=1_{B}$ for all $a\in A$. Then $\overline{\top}$ is the top element in $B^{A}$.

(5) Define a map $\top \colon A\longrightarrow B$ as follows:
\[\top(a)=\left\{
               \begin{array}{ll}
                 0_{B}, & \hbox{$a=0_{A}$} \\
                 1_{L}, & \hbox{otherwise}
               \end{array}
             \right.
\]
for all $a\in A$. Then $\top$ is the top element in $\mathcal{S}(A,B)$. $\underline{\perp}$ is the bottom element in $\mathcal{S}(A,B)$.

(6) Define a map $\perp \colon A\longrightarrow B$ as follows:
\[\top(a)=\left\{
               \begin{array}{ll}
                 1_{B}, & \hbox{$a=1_{A}$} \\
                 0_{L}, & \hbox{otherwise}
               \end{array}
             \right.
\]
for all $a\in A$. Then $\perp$ is the bottom element in $\mathcal{M}(A,B)$. $\overline{\top}$ is the top element in $\mathcal{M}(A,B)$.
\end{remark}

\begin{definition}\label{def2.2}  A \emph{quantale} is a complete lattice  $L$ with an associative binary operation $\ast$ satisfying:
\[a\ast\Big(\bigvee\limits_{i\in I}b_{i}\Big)=\bigvee\limits_{i\in
I}(a\ast b_{i})\  \mbox{and}\ \Big(\bigvee\limits_{i\in
I}b_i\Big)\ast a=\bigvee\limits_{i\in I}(b_i\ast a)\]
for all $ {a}\in{L},\{{b}_{i}\}_{i\in I}\subseteq{L}$.
\end{definition}
A quantale $L$ is said to be \emph{unital} provided that there exists an element $e_{L}\in L$ such that ${a}\ast e_{L}={a}$ and
$e_{L}\ast{a}={a}$ for all $a\in L$. An element $x$ is called a {\it left unit} if $x\ast a= a$ for all $a\in L$. Dually, we can give the definition of {\it right unit}. A nonempty subset $T\subseteq L$ is called a \emph{subquantale} of $L$ if it is closed under all sups and operation in $L$.
 A  $L$ is said to be  {\it simple} if any surjective morphism of  quantales is either an isomorphism or a constant morphism.

\begin{definition}\label{def2.3} Let $M$ and $N$ be quantales. A \emph{quantale homomorphism} $f\colon M \longrightarrow N$ is a semigroup
homomorphism  such that $f$ preserves arbitrary sups, that is,
$f(\bigvee_{i\in I}x_{i})$ $=\bigvee_{i\in I}f(x_{i})$ for all  $\{{x}_{i}\}_{i\in I}\subseteq{M}$.
\end{definition}

\begin{example}\label{def2.4} Let $L$ be a complete lattice. Then $(\mathcal{S}(L),\circ)$ be a quantale.
\end{example}

$(\mathcal{S}(L),\circ)$ is very important in quantale theory. Since $(\mathcal{S}(L),\circ)$ is unital, every quantale can be embedded into a unital quantale. $(\mathcal{S}(L),\circ)$ is widely used as a simple quantale. Besides, $(\mathcal{S}(L),\circ)$ is a girard quantale whenever $L$ is a completely distributive lattice.

\begin{definition}\label{def2.5}
Let $(L,\ast_{L})$ be a quantale.\ A map $j\colon L\longrightarrow L$ is called a \emph{quantale nucleus} on $L$ provided that for all $a,b\in L$,

(1) $a\leq j(a)$;

(2) $a\leq b$ implies $j(a)\leq j(b)$;

(3) $j\circ j(a)=j(a)$;

(4) $j(a)\ast_{L} j(b)\leq j(a\ast_{L} b)$.
\end{definition}
Let $L_{j}$ denote the set $\{l\in L\mid j(l)=l\}$. $L_{j}$ is closed under the arbitrary infs. If $(L,\ast_{L})$ is a quantale, then $(L_{j},\ast_{L_{j}})$ is a quantale where the $\ast_{L_{j}}$ is given by $a\ast_{L_{j}} b=j(a\ast_{L} b)$ for all $a,b\in L_{j}$.

\begin{definition} Let $(L,\ast)$ be a quantale. A subset $T\subseteq L$ is called a quantale quotient of $L$ if there exists a nucleus $j$ on $L$ such that $T=L_{j}$.
\end{definition}

\section{The quantale of order-preserving maps}

From Raney (Definition 3 and Theorem 1 in \cite{Raney1953}), we say that $x$ is {\it wedge-below} to $y$ in a poset $P$, and we write $x\triangleleft y$, for any
subset $A\subseteq P$ with $y\leq \bigvee A$ there is an $a\in A$ such that $x\leq a$.
\begin{definition}\label{def3.1}(\cite{Raney1953}) Let $L$ be a complete lattice. $L$ is called a {\it completely distributive lattice} if $x=\bigvee\{y\mid y\triangleleft x\}$ for all $x$ in $L$.\end{definition}
\begin{remark}\label{rem3.2}Let $L$ be a completely distributive lattice.

(1) $0$ is not wedge-below to $0$ and $0$ is wedge-below to the other elements in $L$.

(2) If $x_{1}\triangleleft y$ and $x_{2}\leq x_{1}$ in $L$, then $x_{2}\triangleleft y$.

(3) If $x\triangleleft y_{1}$ and $y_{1}\leq y_{2}$ in $L$, then $x\triangleleft y_{2}$.

(4) If $x\triangleleft y$, then $x\leq y$.

(5) $x\triangleleft \bigvee_{i\in I}y_{i}$ in $L$ if and only if there exists $i_{0}\in I$ such that $x\triangleleft y_{i_{0}}$.
\end{remark}

Let $A,B,C$ be complete lattices and $f\in B^{A},g\in C^{B}$. The new composition of order-preserving maps $g\cdot f\colon A\longrightarrow C$ is
defined by
\[(g\cdot  f)(a) := \bigvee\{c \in C \mid b \in B: b\triangleleft f(a)\ \mbox{and}\ c\triangleleft g(b)\}\]
for all $a\in A$.
\begin{remark}\label{rem3.3} Let $A, B$ and $C$ be  complete lattices. If $f\in B^{A}$ and $g\in C^{B}$, then $g\cdot f\leq g\circ f$.
\end{remark}
\begin{proposition}\label{pro3.4} Let $A,B$ and $C$ be complete lattices. If $f\in B^{A}$ and $g\in C^{B}$, then $g\cdot  f\in C^{A}$.
\begin{proof} Let $x,y\in A$ with $x\leq y$. Then $f(x)\leq f(y)$, which implies $\{b\triangleleft f(x)\}\subseteq \{b\triangleleft f(y)\}$.\ It follows that $\{c\triangleleft g(b)\mid b \triangleleft f(x)\}\subseteq \{c\triangleleft g(b)\mid b \triangleleft f(y)\}$ and hence
$(g\cdot  f)(x)\leq (g\cdot  f)(y)$.\end{proof}
\end{proposition}
\begin{remark}\label{rem3.5} (1) Let $A, B$ be  complete lattices and $C$ a completely distributive lattices.
If $f\in B^{A}$ and $g\in C^{B}$, then $(g\cdot  f)(a)=\bigvee\{g(b)\mid b \in B:b\triangleleft f(a)\}$ for all $a\in A$.

(2) Let $A, B$ and $C$ be complete lattices. If $f\in B^{A}$ and $g\in C^{B}$, then $g\cdot f\leq g\circ f$.
\end{remark}

Note that all mentioned complete lattices in the following paper are always completely distributive lattices.

\begin{proposition}\label{pro3.6} Let $A,B$ and $C$ be complete lattices.

(1) If $f\in \mathcal{S}(A,B)$ and $g\in C^{B}$, then $g\cdot  f\in \mathcal{S}(A,C)$.

(2) If $f\in \mathcal{S}(A,B)$ and $g\in \mathcal{S}(B,C)$, then $g\cdot  f=g\circ f$.
\begin{proof} (1) Let $\{a_{i}\}\subseteq A$. Then  $(g\cdot f)(\bigvee_{i\in I} a_{i})=\bigvee\{g(b)\mid b \in B:b\triangleleft f(\bigvee_{i\in I} a_{i})\}$. Since $f\in \mathcal{S}(A,B)$, $(g\cdot f)(\bigvee_{i\in I} a_{i})=\bigvee\{g(b)\mid b \in B:b\triangleleft \bigvee_{i\in I}f( a_{i})\}=\bigvee\{g(b)\mid b \in B: \exists i_{0}\in I, b\triangleleft f( a_{i_{0}})\}=\bigvee_{i\in I}\bigvee\{g(b)\mid b \in B:b\triangleleft f( a_{i})\}$.

(2) By(1), $g\cdot  f\in \mathcal{S}(A,C)$. For all $a\in A$,  since $g\in \mathcal{S}(B,C)$, we have that $(g\cdot  f)(a)=\bigvee\{g(b)\mid b \in B:b\triangleleft f(a)\}=g(\bigvee\{b \in B:b\triangleleft f(a)\})=g(f(a))$. Thus $g\cdot  f=g\circ f$.\end{proof}\end{proposition}

\begin{proposition}\label{pro3.7} Let $A,B$ and $C$ be complete lattices. If $f\in B^{A}$, $g\in C^{B}$ and $h\in D^{C}$, then  $(h\cdot g)\cdot  f = h\cdot (g \cdot f)$.
\begin{proof}Firstly let us show that $((h\cdot g)\cdot  f)(a)\leq  (h\cdot (g \cdot f))(a)$ for all $a\in A$. For each $x\in D$ with $x\triangleleft ((h\cdot g)\cdot  f)(a)$,  we have $x\triangleleft \bigvee\{d\mid  b \in B:b\triangleleft f(a), d\triangleleft (h\cdot g)(b)\}$.
Then there exist $d_{0}$ and $b_{0}$ such that $x\triangleleft d_{0}$, $b_{0}\triangleleft f(a)$  and $d_{0}\triangleleft (h\cdot g)(b_{0})$. $(h\cdot g)(b_{0})=\bigvee\{d\mid  c \in C:c\triangleleft g(b_{0}), d\triangleleft h(c)\}$.
Hence, there exist $d_{1}\in D$, $c_{0}\in C$ such that $d_{0}\triangleleft d_{1}$, $c_{0}\triangleleft g(b_{0})$ and $d_{1}\triangleleft h(c_{0})$. $(h\cdot (g \cdot f))(a)=\bigvee\{d\mid  c \in C:c\triangleleft (g\cdot f)(a), d\triangleleft h(c)\}$. $(g\cdot f)(a)=\bigvee\{c\mid b \in B:b\triangleleft f(a), c\triangleleft g(b)\}$. By the above discussion, $c_{0}\triangleleft (g\cdot f)(a)$, $d_{1}\triangleleft h(c_{0})$ and $d_{1}\leq (h\cdot (g \cdot f))(a)$ which implies that $x\leq (h\cdot (g \cdot f))(a)$. Thus  $(h\cdot g)\cdot  f\leq h\cdot (g \cdot f)$.

Next let us prove the inequality $((h\cdot g)\cdot  f)(a)\geq (h\cdot (g \cdot f))(a)$ for all $a\in A$. For each $x\triangleleft (h\cdot (g \cdot f))(a)$,  we have $x\triangleleft \bigvee\{d\mid  c \in C:c\triangleleft (g \cdot f)(a), d\triangleleft h(c)\}$. Then there exist $d_{0}$ and $c_{0}$ such that $x\triangleleft d_{0}$, $c_{0}\triangleleft (g \cdot f)(a)$ and $d_{0}\triangleleft h(c_{0})$. Since $(g\cdot f)(a)=\bigvee\{c\mid b \in B:b\triangleleft f(a), c\triangleleft g(b)\}$, we have that there exist $c_{1}\in C$, $b_{0}\in B$ such that $c_{0}\triangleleft c_{1}$, $b_{0}\triangleleft f(a)$ and $c_{1}\triangleleft g(b_{0})$. So $c_{0}\triangleleft h(b_{0})$, $b_{0}\triangleleft f(a)$ and $d_{0}\triangleleft h(c_{0})$.
Since $((h\cdot g)\cdot  f)(a)=\bigvee\{d\mid  b \in B:b\triangleleft f(a), d\triangleleft (h\cdot g)(b)\}$ and $(h\cdot g)(b)=\bigvee\{d\mid  c \in C:c\triangleleft g(b), d\triangleleft h(c)\}$, which implies $x\leq ((h\cdot g)\cdot  f)(a)$. Thus $(h\cdot g)\cdot  f\geq h\cdot (g \cdot f)$.
\end{proof}\end{proposition}

\begin{theorem}\label{th3.8} Let $L$ be a complete lattices. Then $(L^{L},\cdot)$ is a quantale.
\begin{proof} For all $\{f_{i}\}_{i\in I}\subseteq L^{L}$ and $g\in L^{L}$. For any $a\in L$, $(g\cdot (\bigvee_{i\in I}f_{i}))(a)=\bigvee\{g(b)\mid b \in L: b\triangleleft (\bigvee_{i\in I}f_{i})(a)\}=\bigvee\{g(b)\mid b\in L: b\triangleleft \bigvee_{i\in I}(f_{i}(a))\}=\bigvee\{g(b)\mid b \in L: \exists i_{0}\in I,\ b\triangleleft f_{i_{0}}(a)\}=\bigvee_{i\in I}\bigvee\{g(b)\mid b\in L:  b\triangleleft f_{i}(a)\}=(\bigvee_{i\in I}(g\cdot f_{i}))(a)$. Moreover, $( (\bigvee_{i\in I}f_{i})\cdot g)(a)=\bigvee\{(\bigvee_{i\in I}f_{i})(b)\mid b \in L: b\triangleleft g(a)\}=\bigvee\{\bigvee_{i\in I}(f_{i}(b))\mid b \in L: b\triangleleft g(a)\}=\bigvee_{i\in I}\bigvee\{f_{i}(b)\mid b \in L: b\triangleleft g(a)\}=(\bigvee_{i\in I}(f_{i}\cdot g))(a)$. Thus, $(L^{L},\cdot)$ is a quantale.
\end{proof}
\end{theorem}
By Proposition \ref{pro3.6}(2), we can immediately obtain the following corollary.

\begin{corollary}\label{cor3.9} $(\mathcal{S}(L),\circ)$ is a sub-quantale of $(L^{L},\cdot)$.
\end{corollary}

By the definition of $(L^{L},\cdot)$, one can easily check the following properties.

\begin{corollary}\label{cor3.10} (1) The identity map is the left unit of $(L^{L},\cdot)$.

(2) $\overline{\top}\cdot f=\overline{\top}$ for all $f\in \mathcal{M}(L)$ with $f(0_{A})\neq 0_{B}$.

(3) $\overline{\top}\cdot \top=\top$.

(4) $\top\cdot \overline{\top}=\overline{\top}$.

(5) Define a map $f_{a}\colon L\longrightarrow L$ by $f_{a}(0)=0$ and $f_{a}(x)=a$ for all $x\in L\setminus\{0\}$. Then $f_{a}\cdot f_{a}\neq f_{a}$.
\end{corollary}
\begin{example}\label{exam3.11} The identity map is not the right unit of $(L^{L},\cdot)$. Clearly, the identity map $id_{L}\in \mathcal{S}(L)$, then $\overline{\top}\cdot id_{L}=\top$.
\end{example}

By the above corollary and example, we know that  $(L^{L},\cdot)$ is not unital, commutative, idempotent and two-sided. So $(L^{L},\cdot)$ is a rather special quantale which can be used to put counterexample in some cases.

\begin{proposition}\label{pro3.12} Let $A,B$ and $C$ be complete lattices. For  each $f\in  B^{A}$, define a map $\psi(f)\colon A \longrightarrow  B$ by
\[\psi(f)(a) := \bigvee\{b\in B\mid s\triangleleft a\ \mbox{and}\ b\triangleleft f(s)\},\]
for all $a\in A$. Then $\psi(f)\in \mathcal{S}(A,B)$.
\begin{proof} For any subset $X\subseteq A$, $\psi(f)(\bigvee X)= \bigvee\{b\in B\mid s\triangleleft \bigvee X\ \mbox{and}\ b\triangleleft f(s)\}=\bigvee\{f(s)\mid s\triangleleft \bigvee X\}=\bigvee_{x\in X}\bigvee\{f(s)\mid s\triangleleft x\}=\bigvee_{x\in X}\psi(f)(x)$. Thus, we have that $\psi(f)\in \mathcal{S}(A,B)$.\end{proof}\end{proposition}
By the above proposition, one can easily check the following properties.
\begin{corollary}\label{cor3.13} (1) For  each $f\in  B^{A}$, $\psi(f)(a)=\bigvee\{f(s)\mid s\triangleleft a\}$ for all $a\in A$.

(2) For  each $f\in  B^{A}$, $\psi(f)\leq f$.

(3) For  each $f\in  \mathcal{S}(A,B)$, $\psi(f)=f$.

(4) For  each $f,g\in  B^{A}$, $\psi(f)\leq g$ if and only if $f\leq g$.

\end{corollary}

\begin{proposition}\label{pro3.14} Define a map $k\colon L^{L}\longrightarrow \mathcal{S}(L)$ by $k(f)=\psi(f)$. Then $k$ is a surjective quantale homomorphism.
\begin{proof} By Proposition \ref{pro3.12}, $k$ is well-defined. Let $g,f\in L^{L}$, for all $a\in L$, on the one hand, $\psi (g\cdot f)(a)=\bigvee\{(g\cdot f)(s)\mid s\triangleleft a\}=\bigvee\{g(t)\mid t\triangleleft f(s), s\triangleleft a\}$, on the other hand, $(\psi (g)\circ \psi(f))(a)=\psi (g)(\psi(f)(a))=\bigvee\{g(s)\mid s\triangleleft \psi(f)(a)\}=\bigvee\{g(s)\mid s\triangleleft \bigvee\{f(t)\mid t\triangleleft a\}\}=\bigvee\{g(s)\mid s\triangleleft \bigvee\{f(t)\mid t\triangleleft a\}\}=\bigvee\{g(s)\mid s\triangleleft f(t), t\triangleleft a\}$. So $\psi (g\cdot f)=\psi (g)\circ \psi (f)$.

For all $\{f_{i}\}_{i\in I}\subseteq L^{L}$, $k(\bigvee_{i\in I}f_{i})(a)=\psi(\bigvee_{i\in I}f_{i})(a)=\bigvee\{(\bigvee_{i\in I}f_{i})(s)\mid s\triangleleft a\}=\bigvee\{\bigvee_{i\in I}(f_{i}(s))\mid s\triangleleft a\}=\bigvee_{i\in I}\bigvee\{f_{i}(s)\mid s\triangleleft a\}=\bigvee_{i\in I}(\psi(f_{i})(a))=(\bigvee_{i\in I} \psi(f_{i}))(a)$. By the above discussion, $k$ is a quantale homomorphism. Moreover, it is easy to see that $\psi$ is surjective by Corollary  \ref{cor3.13}(3).\end{proof}\end{proposition}

By Corollary \ref{cor3.9}, Proposition \ref{pro3.14}, we can obtain the following corollary.
\begin{corollary}\label{cor3.15} (1) $(\mathcal{S}(L),\circ)$ is a retract of $(L^{L},\cdot)$.

(2) If $h\colon L^{L}\longrightarrow M$ is a quantale homomorphism such that $h(\mathcal{S}(L))=M$, then $M$ is quantale isomorphic to $\mathcal{S}(L)$.
\end{corollary}

In the following, we give a representation on $(\mathcal{S}(L),\circ)$.

\begin{theorem}\label{th3.16} There exists a nucleus $j$ on $(L^{L},\cdot)$ such that $(\mathcal{S}(L),\circ)\cong ((L^{L})_{j},\cdot_{j})$.
\begin{proof} Define a map $k^{\ast}\colon \mathcal{S}(L)\longrightarrow L^{L}$ be $k^{\ast}(f)=\bigvee\{g\in L^{L}\mid k(g)\leq f\}$. $k^{\ast}$ is the right adjoint of $k$. Then $id_{L^{L}} \leq k^{\ast}\circ k$ and $id_{\mathcal{S}(L)}=k\circ k^{\ast}$. Let $j:\equiv k^{\ast}\circ k$. Clearly, $j$ is an order-preserving map. Next, we show that $j$ is a nucleus on $L^{L}$. Firstly, since $id_{L^{L}} \leq k^{\ast}\circ k$, we have that $f\leq j(f)$ for all $f\in L^{L}$. Secondly, $j(j(f))=(k^{\ast}\circ k)(k^{\ast}\circ k(f))=(k^{\ast}\circ k\circ k^{\ast}\circ k)(f)=k^{\ast}\circ k(f)=j(f)$. Finally, for any $f,g\in L^{L}$, we have that $j(g\cdot f)=(k^{\ast}\circ k)(g\cdot f)=k^{\ast}(k(g)\circ k(f))=\bigvee\{h\in L^{L}\mid k(h)\leq k(g)\circ k(f)\}\geq \bigvee\{s\in L^{L}\mid k(s)\leq k(g)\}\cdot \bigvee\{t\in L^{L}\mid k(t)\leq k(f)\}=j(g)\cdot j(f)$. Thus, $j$ is a nucleus on $L^{L}$ and  $((L^{L})_{j},\cdot_{j})$ is a quantale. Define a map $\phi\colon \mathcal{S}(L)\longrightarrow  (L^{L})_{j}$ by $\phi(f)=k^{\ast}(f)$ for all $f\in \mathcal{S}(L)$. $j(\phi(f))=(k^{\ast}\circ k)(k^{\ast}(f))=k^{\ast}(f)$. For all $f\in (L^{L})_{j}$, there exists $k(f)$ in $\mathcal{S}(L)$ such that $\phi(k(f))=(k^{\ast}\circ k)(f)=j(f)=f$. So $\phi$ is an ordered isomorphism. For any $g,f\in \mathcal{S}(L)$, $\phi(g\circ f)=k^{\ast}(k(g)\circ k(f))=j(j(g)\cdot j(f))=\phi(g)\cdot_{j}\phi(f)$. Thus, $(\mathcal{S}(L),\circ)\cong ((L^{L})_{j},\cdot_{j})$. \end{proof}\end{theorem}

\section{The co-quantale of order-preserving maps}

\begin{definition}\label{def4.1}  A \emph{co-quantale} is a complete lattice  $L$ with an associative binary operation $\ast$ satisfying:
\[a\ast\Big(\bigwedge\limits_{i\in I}b_{i}\Big)=\bigwedge\limits_{i\in
I}(a\ast b_{i})\  \mbox{and}\ \Big(\bigwedge\limits_{i\in
I}b_i\Big)\ast a=\bigwedge\limits_{i\in I}(b_i\ast a)\]
for all $ {a}\in{L},\{{b}_{i}\}_{i\in I}\subseteq{L}$.
\end{definition}

The co-quantales is a contravariant notion of the quantales. In \cite{RZ2021}, co-quantale plays an important role in the study of valued logic.

\begin{proposition}\label{pro4.2} Let $L$ be a complete lattice. Then $(\mathcal{M}(L),\circ)$ is a co-quantale.
\begin{proof} We know that $\circ$ is an associative binary operator. Next we show that $(\bigwedge_{i\in I}g_{i})\circ f=\bigwedge_{i\in I}(g_{i}\circ f)$. For all $x\in L$, we have that $((\bigwedge_{i\in I}g_{i})\circ f)(x)=(\bigwedge_{i\in I}g_{i})(f(x))=\bigwedge_{i\in I}(g_{i}(f(x)))=\bigwedge_{i\in I}(g_{i}\circ f)(x)=(\bigwedge_{i\in I}(g_{i}\circ f))(x)$. Similarly, we can prove  $f\circ(\bigwedge_{i\in I}g_{i})=\bigwedge_{i\in I}(f\circ g_{i})$. Thus, $(\mathcal{M}(L),\circ)$ is a co-quantale.\end{proof}
\end{proposition}

In this section, we aim to construct a new composition on $L^{L}$, which allows us to endow $L^{L}$ a quantale structure.

By \cite{GKLMS2003}, we know that $L$ is completely distributive lattice if and only if $L^{op}$ is completely distributive.
We say $x$ is {\it co-wedge-below} to $y$ in $L$, and we write $x\triangleleft^{co} y$, for any
subset $A\subseteq L$ with $\bigwedge A\leq y$ there is an element $a\in A$ such that $a\leq x$.
\begin{remark}\label{rem4.3} (1) If $x\triangleleft^{co}y$ in $L$, then $y\leq x$.

(2) If $x_{1}\triangleleft^{co}y$ and $x_{1}\leq x_{2}$, then $x_{2}\triangleleft^{co}y$.

(3) If $x\triangleleft^{co}y_{1}$ and $y_{2}\leq y_{1}$, then $x\triangleleft^{co}y_{2}$.

(4)  $x\triangleleft^{co}\bigwedge_{i\in I}y_{i}$ if and only if there exists $i_{0}\in I$ such that $x\triangleleft^{co} y_{i_{0}}$.

(5) $L$ is  a  completely distributive lattice if $x=\bigwedge\{y\mid y\triangleleft^{co} x\}$ for all $x$ in $L$.
\end{remark}

Let  $A,B,C$ be complete lattices and $f\in B^{A},g\in C^{B}$. Another composition of order-preserving maps $g\bullet f\colon A\longrightarrow C$ is
defined by
\[(g\bullet  f)(a) := \bigwedge\{c \in C \mid b \in B: b\triangleleft^{co} f(a)\ \mbox{and}\ c\triangleleft^{co} g(b)\}\]
for all $a\in A$.

\begin{proposition}\label{pro4.4} If $f\in B^{A}$ and $g\in C^{B}$, then $g\bullet  f\in C^{A}$.
\begin{proof} Let $x,y\in A$ with $x\leq y$. Then $f(x)\leq f(y)$, which implies $\{b\triangleleft^{co} f(y)\}\subseteq \{b\triangleleft^{co} f(x)\}$.\ It follows that $\{c\triangleleft^{co} g(b)\mid b \triangleleft^{co} f(y)\}\subseteq \{c\triangleleft^{co} g(b)\mid b \triangleleft^{co} f(x)\}$ and
$\bigwedge\{c\triangleleft^{co} g(b)\mid b \triangleleft^{co} f(y)\}\geq \bigwedge\{c\triangleleft^{co} g(b)\mid b \triangleleft^{co} f(x)\}$.
Thus $(g\bullet  f)(x)\leq (g\bullet  f)(y)$.
\end{proof}
\end{proposition}
\begin{remark}\label{rem4.5} (1) If $f\in B^{A}$ and $g\in C^{B}$, then $(g\bullet  f)(a)=\bigwedge\{g(b)\mid b \in B:b\triangleleft^{co} f(a)\}$ for all $a\in A$.

(2) If $f\in B^{A}$ and $g\in C^{B}$, then $g\circ f\leq g\bullet f$.
\end{remark}

\begin{proposition}\label{pro4.6} (1) If $f\in \mathcal{M}(A,B)$ and $g\in C^{B}$, then $g\bullet  f\in \mathcal{M}(A,C)$.

(2) If $f\in \mathcal{M}(A,B)$ and $g\in \mathcal{M}(B,C)$, then $g\bullet  f=g\circ f$.

(3) If $f\in B^{A}$ and $g\in C^{B}$, then $g\cdot  f\leq g\bullet  f$.

(4) If $f\in \mathcal{M}(A,B)\bigcap\mathcal{S}(A,B)$, $g\in \mathcal{M}(B,C)\bigcap\mathcal{S}(B,C)$, then $g\bullet  f=g\cdot  f$.
\begin{proof} (1) Let $\{a_{i}\}\subseteq A$. Then  $(g\bullet f)(\bigwedge_{i\in I} a_{i})=\bigwedge\{g(b)\mid b \in B:b\triangleleft^{co} f(\bigwedge_{i\in I} a_{i})\}$. Since $f\in \mathcal{S}(A,B)$, $(g\bullet f)(\bigwedge_{i\in I} a_{i})=\bigwedge\{g(b)\mid b \in B:b\triangleleft^{co} \bigvee_{i\in I}f( a_{i})\}=\bigvee\{g(b)\mid b \in B: \exists i_{0}\in I, b\triangleleft^{co} f( a_{i_{0}})\}=\bigwedge_{i\in I}\bigwedge\{g(b)\mid b \in B:b\triangleleft^{co} f( a_{i})\}$.

(2) By(1), $g\bullet f\in \mathcal{S}(A,C)$. For all $a\in A$,  since $g\in \mathcal{S}(B,C)$, $(g\cdot  f)(a)=\bigwedge\{g(b)\mid b \in B:b\triangleleft^{co} f(a)\}=g(\bigwedge\{b \in B:b\triangleleft^{co} f(a)\})=g(f(a))$. Thus $g\bullet  f=g\circ f$.

(3) For any $a\in A$, $(g\cdot f)(a)=\bigvee\{g(b)\mid b\triangleleft f(a)\}\leq\bigwedge\{g(x)\mid x\triangleleft^{co}f(a)\}=(g\bullet f)(a)$.

(4) By(2) and Proposition \ref{pro3.6}(2), it is easy to prove it.\end{proof}
\end{proposition}

\begin{proposition}\label{pro4.7} Let  $A,B,C$ and $D$ be complete lattices. If $f\in B^{A}$, $g\in C^{B}$ and $h\in D^{C}$, then  $(h\bullet g)\bullet  f = h\bullet (g \bullet f)$.
\begin{proof}Firstly let us show that $((h\bullet g)\bullet  f)(a)\geq  (h\bullet (g \bullet f))(a)$ for all $a\in A$. For each $x\in D$ with $x\triangleleft^{co} ((h\bullet g) \bullet f)(a)$, we have that  $x\triangleleft^{co} \bigwedge\{d\mid  b \in B: b\triangleleft^{co} f(a),\ d\triangleleft^{co} (h\bullet g)(b)\}$. Then there exist $d_{0}$ and $b_{0}$ such that $x\triangleleft^{co} d_{0}$, $b_{0}\triangleleft^{co} f(a)$  and $d_{0}\triangleleft^{co} (h\bullet g)(b_{0})$.  $(h\bullet g)(b_{0})=\bigwedge\{d\mid  c \in C:c\triangleleft^{co} g(b_{0}), d\triangleleft^{co} h(c)\}$.
Hence, there exist $d_{1}\in D$, $c_{0}\in C$ such that $d_{0}\triangleleft^{co} d_{1}$, $c_{0}\triangleleft^{co} g(b_{0})$ and $d_{1}\triangleleft^{co} h(c_{0})$. $(h\bullet (g \bullet f))(a)=\bigwedge\{d\mid  c \in C:c\triangleleft^{co} (g\cdot f)(a),\ d\triangleleft^{co} h(c)\}$. $(g\bullet f)(a)=\bigwedge\{c\mid b \in B:b\triangleleft^{co} f(a),\ c\triangleleft^{co} g(b)\}$. By the above discussion, $c_{0}\triangleleft^{co} (g\bullet f)(a)$ and $d_{1}\triangleleft^{co} h(c_{0})$, $d_{1}\geq (h\bullet (g \bullet f))(a)$ which implies that $x\geq (h\bullet (g \bullet f))(a)$. Thus  $(h\bullet g)\bullet f\geq h\bullet (g \bullet f)$.

Next let us prove the inequality $((h\bullet g)\bullet  f)(a)\leq (h\bullet (g \bullet f))(a)$ for all $a\in I$. For each $x\in D$ with $x\triangleleft^{co} (h\bullet (g \bullet f))(a)$, we have $x\triangleleft^{co} \bigwedge\{d\mid  c \in C:c\triangleleft^{co} (g \bullet f)(a),\ d\triangleleft^{co} h(c)\}$. Then there exist $d_{0}$ and $c_{0}$ such that $x\triangleleft^{co} d_{0}$, $c_{0}\triangleleft^{co} (g \bullet f)(a)$ and $d_{0}\triangleleft^{co} h(c_{0})$. Since $(g\bullet f)(a)=\bigwedge\{c\mid b \in B:b\triangleleft^{co} f(a),\ c\triangleleft^{co} g(b)\}$, we have that there exist $c_{1}\in C$, $b_{0}\in B$ such that $c_{0}\triangleleft^{co} c_{1}$, $b_{0}\triangleleft^{co} f(a)$ and $c_{1}\triangleleft^{co} g(b_{0})$. So $c_{0}\triangleleft^{co} h(b_{0})$, $b_{0}\triangleleft^{co} f(a)$ and $d_{0}\triangleleft^{co} h(c_{0})$.
Since $((h\bullet g)\bullet  f)(a)=\bigwedge\{d\mid  b \in B:b\triangleleft^{co} f(a), d\triangleleft^{co} (h\bullet g)(b)\}$ and $(h\bullet g)(b)=\bigwedge\{d\mid  c \in C:c\triangleleft^{co} g(b), d\triangleleft^{co} h(c)\}$, which implies $x\geq ((h\bullet g)\bullet  f)(a)$. Thus $(h\bullet g)\bullet  f\leq h\bullet (g \bullet f)$.\end{proof}
\end{proposition}

\begin{theorem}\label{th4.8} $(L^{L},\bullet)$ is a co-quantale.
\begin{proof} For all $\{f_{i}\}_{i\in I}\subseteq L^{L}$ and $g\in L^{L}$. For any $a\in L$, $(g\bullet (\bigwedge_{i\in I}f_{i}))(a)=\bigwedge\{g(b)\mid b \in L: b\triangleleft^{co} (\bigwedge_{i\in I}f_{i})(a)\}=\bigwedge\{g(b)\mid b\in L: b\triangleleft^{co} \bigwedge_{i\in I}(f_{i}(a))\}=\bigwedge\{g(b)\mid b \in L: \exists i_{0}\in I,\ b\triangleleft^{co} f_{i_{0}}(a)\}=\bigwedge_{i\in I}\bigwedge\{g(b)\mid b\in L:  b\triangleleft^{co} f_{i}(a)\}=(\bigwedge_{i\in I}(g\bullet f_{i}))(a)$. Moreover, $( (\bigwedge_{i\in I}f_{i})\bullet g)(a)=\bigwedge\{(\bigwedge_{i\in I}f_{i})(b)\mid b \in L: b\triangleleft^{co} g(a)\}=\bigwedge\{\bigwedge_{i\in I}(f_{i}(b))\mid b \in L: b\triangleleft^{co} g(a)\}=\bigwedge_{i\in I}\bigwedge\{f_{i}(b)\mid b \in L: b\triangleleft^{co} g(a)\}=(\bigwedge_{i\in I}(f_{i}\bullet g))(a)$. Thus, $(L^{L},\bullet)$ is a co-quantale.
\end{proof}
\end{theorem}

\begin{proposition}\label{pro4.9} For  each $f\in  B^{A}$, define $\varphi(f)\colon A \longrightarrow  B$ by
\[\varphi(f)(a) := \bigwedge\{b\in B \mid s\triangleleft^{co} a\ \mbox{and}\ b\triangleleft^{co} f(s)\},\]
then $\varphi(f)\in \mathcal{M}(A,B)$.
\begin{proof} For any subset $X\subseteq A$, $\varphi(f)(\bigwedge X)= \bigwedge\{b\in B \mid s\triangleleft^{co}\bigwedge X\ \mbox{and}\ b\triangleleft^{co} f(s)\}=\bigwedge\{f(s) \mid s\triangleleft^{co}\bigwedge X\}=\bigwedge_{x\in X}\bigwedge\{f(s)\mid s\triangleleft^{co} x\}=\bigwedge_{x\in X}\varphi(f)(x)$. Thus, $\varphi(f)\in \mathcal{M}(A,B)$.\end{proof}
\end{proposition}

\begin{corollary}\label{cor4.10} (1) For  each $f\in  B^{A}$, $\varphi(f)(a)=\bigwedge\{f(s)\mid s\triangleleft^{co} a\}$ for all $a\in A$.

(2) For  each $f\in  B^{A}$, $f\leq \varphi(f)$.

(3) For  each $f\in  \mathcal{M}(A,B)$, $\varphi(f)=f$.
\end{corollary}

\begin{definition}\label{def4.11}
Let $(L,\ast_{L})$ be a co-quantale.\ A map $j\colon L\longrightarrow L$ is called a \emph{co-quantale nucleus} on $L$ provided that for all $a,b\in L$,

(1) $a\leq j(a)$;

(2) $a\leq b$ implies $j(a)\leq j(b)$;

(3) $j\circ j(a)=j(a)$;

(4) $j(a)\ast_{L} j(b)\leq j(a\ast_{L} b)$.
\end{definition}

\begin{proposition}\label{pro4.12} Let $(L,\ast)$ be a co-quantale and $j$  a nucleus on $L$. Then $(L_{j},\ast_{j})$ is a co-quantale.
\begin{proof} Since $L_{j}$ is closed under  the arbitrary infs, it is a complete lattice. For any $x,y$ and $z$ in $L$, we have that
$(x\ast_{j}y)\ast_{j}z=j(j(x\ast y)\ast z)\geq j((x\ast y)\ast z)$. Besides, $j((x\ast y)\ast z)\leq j(j(x\ast y)\ast j(z))\leq j(j((x\ast y)\ast z))=j((x\ast y)\ast z)$. Hence, $(x\ast_{j}y)\ast_{j}z=j((x\ast y)\ast z)$. $x\ast_{j}(y\ast_{j} z)=j(x\ast j(y\ast z))\leq j(j(x)\ast j(y\ast z))\leq
j(j(x\ast (y\ast z)))=j(x\ast (y\ast z))$. $j(x\ast (y\ast z))\leq j(x\ast j(y\ast z))$. Hence, $x\ast_{j}(y\ast_{j} z)=j(x\ast j(y\ast z))$. Since $\ast$ is an associative binary operator, $(x\ast_{j}y)\ast_{j}z=x\ast_{j}(y\ast_{j} z)$. For all $\{x_{i}\}_{i\in I}\subseteq L$ and $x\in L$, we have that $(\bigwedge_{i\in I} x_{i})\ast_{j} x=j((\bigwedge_{i\in I} x_{i})\ast x)=j(\bigwedge_{i\in I}( x_{i}\ast x))=\bigwedge_{i\in I} j(x_{i}\ast x)=\bigwedge_{i\in I} (x_{i}\ast_{j} x)$. Similarly, we can prove that $x\ast_{j}(\bigwedge_{i\in I} x_{i})=\bigwedge_{i\in I} (x\ast_{j}x_{i})$.
Thus $(L_{j},\ast_{j})$ is a co-quantale.
\end{proof}
\end{proposition}

\begin{definition}\label{def4.13} Let $(L,\ast)$ be a co-quantale. A subset $T\subseteq L$ is called a co-quantale quotient of $L$ if there exists a nucleus $j$ on $L$ such that $T=L_{j}$.
\end{definition}

\begin{theorem}\label{th4.14} $(\mathcal{M}(L),\circ)$ is a co-quantale quotient of $(L^{L},\bullet)$.
\begin{proof} Define a map $j\colon L^{L}\longrightarrow L^{L}$ by $j(f)=\varphi(f)$ for all $f\in L^{L}$. Clearly, $j$ is an order-preserving map. By Corollary \ref{cor4.10}(1), $f\leq j(f)$ for all $f\in L^{L}$. $j(j(f))=\varphi(\varphi(f))=\varphi(f)=j(f)$ by Corollary \ref{cor4.10}(3). For any $f,g\in L^{L}$, we have that $j(g\bullet f)(x)=\varphi(g\bullet f)(x)=\bigwedge\{(g\bullet f)(s)\mid s\triangleleft^{co} x\}=\bigwedge_{s\triangleleft^{co} x}\bigwedge\{g(t)\mid t\triangleleft^{co} f(s)\}$. $(j(g)\bullet j(f))(x)=\bigwedge\{j(g)(y)\mid y\triangleleft^{co} j(f)(x)\}=\bigwedge\{j(g)(y)\mid y\triangleleft^{co} \bigwedge\{f(s)\mid s\triangleleft^{co} x\}\}=\bigwedge_{s\triangleleft^{co} x} j(g)(f(s))=\bigwedge_{s\triangleleft^{co} x}\{g(t)\mid t\triangleleft^{co} f(s)\}$. Hence, $j(g)\bullet j(f)=j(g\bullet f)$ and $j$ is a nucleus on $L^{L}$. By Corollary \ref{cor4.10}, $\mathcal{M}(L)=(L^{L})_{j}$. By Proposition \ref{pro4.6}(2), we know that $\bullet_{j}=\circ$. Thus $(\mathcal{M}(L),\circ)$ is a co-quantale quotient of $(L^{L},\bullet)$.
\end{proof}
\end{theorem}

\bibliographystyle{amsplain}

\end{document}